\documentclass[11pt]{amsart}
\usepackage{amssymb,amsfonts,times,verbatim}
\usepackage{hyperref}

\newcommand{\Qu}{{\mathbb H}}
\newcommand{\Co}{{\mathbb C}}
\renewcommand{\Re}{{\mathbb R}}

\renewcommand{\phi}{\varphi}

\newcommand{\grad}{\operatorname{grad}}

\newcommand{\Spin}{\operatorname{Spin}}
\newcommand{\repart}{\operatorname{Re}}
\newcommand{\End}{\operatorname{End}}
\newcommand{\la}{\langle}       
\newcommand{\ra}{\rangle}       

\newtheorem{thm}{Theorem}[section]
\newtheorem{lemma}[thm]{Lemma}
\newtheorem{prop}[thm]{Proposition}

\newtheorem{remark}[thm]{Remark}
\newtheorem{remarks}[thm]{Remarks}
\newtheorem{definition}[thm]{Definition}
\newtheorem{notation}[thm]{Notation}

\newcommand{\Remark}[1]{\begin{remark}{\rm #1}\end{remark}}

\begin{document}

\title
[Dirac eigenvalues for generic metrics on three-manifolds]
{Dirac eigenvalues for generic metrics on three-manifolds}

\author{Mattias Dahl}

\address{
Institutionen f\"or Matematik
\\
Kungliga Tekniska H\"ogskolan 
\\
100 44 Stockholm
\\
Sweden
}
 
\email{dahl@math.kth.se}



\thanks{
The author was supported by the Swedish Foundation for International 
Cooperation in Research and Higher Education ``STINT'' and by the 
European Commission Network ``Geometric Analysis''. 
}            


\begin{abstract}
We show that for generic Riemannian metrics on a closed spin manifold 
of dimension three the Dirac operator has only simple eigenvalues.
\end{abstract}

\maketitle

\section{Introduction}

The spectrum of the Dirac operator has been explicitly computed for quite 
a few compact spin manifolds, for example spherical space forms, flat 
manifolds, spheres with Berger metrics, see 
\cite{hitchin74a, friedrich84a, baer96a, baer96b, pfaeffle00a}. 
These examples exhibit high multiplicities of the eigenvalues, coming 
from the high degree of symmetry of the spaces, a phenomenon one would
not expect for general Riemannian manifolds. The purpose of this note 
is to prove the following theorem.
\begin{thm}
For a generic metric on a compact spin three-manifold the Dirac operator 
has only simple eigenvalues.
\end{thm}

The Dirac operator $D$ (or $D_g$ to show its dependence on the Riemannian 
metric) is a first order elliptic operator acting on sections of the 
spinor bundle $\Sigma M = \Spin(M,g) \times_{\rho} \Sigma$. It is formally 
self-adjoint and on a compact manifold it has discrete real spectrum. 
Here $\Spin(M,g) \rightarrow SO(M,g) $ is a spin structure on $M$ and 
$\rho: \Spin(n) \rightarrow \End(\Sigma)$, $\Sigma = \Co^{2^{[n/2]}}$, 
is the spinor representation. In dimensions $n \equiv 3,4,5 \mod 8$ the 
spinor representation is quaternionic; there is an $\Re$-linear
endomorphism $J$ of $\Sigma$ with $J^2 = -1$, $Ji = -iJ$, which commutes
with the action of $\Spin(n)$. In particular in three dimensions the spinor 
bundle is a quaternionic line bundle and the eigenspaces $E_{\lambda}$ of 
$D$ are quaternionic vectorspaces. So in three dimensions an eigenvalue 
$\lambda$ is simple if the quaternionic dimension 
$\dim_{\Qu} E_{\lambda} = 1$. For the details of the constructions of spin 
geometry here omitted see \cite{lawson-michelsohn89a}. 

For a compact spin manifold $M$ denote by $R(M)$ the space of Riemannian 
metrics on $M$ equipped with the $C^1$-topology. Denote by $S(M)$ the subset 
consisting of Riemannian metrics for which all eigenvalues of the Dirac 
operator are simple. Recall that a subset of a topological space is called 
residual if it contains a countable intersection of open and dense sets. The 
precise theorem we will prove is the following.
\begin{thm} \label{Mainthm}
Let $M$ be a compact spin manifold of dimension three. Then $S(M)$ is 
residual.
\end{thm}

It is reasonable to conjecture that the conclusion of Theorem \ref{Mainthm} 
holds in any dimension. The restriction to dimension three comes from using 
the simple approach of making conformal variations of the metric. Perhaps 
the techniques of \cite{baer-dahl01ppa, maier97a} can be used to extend the 
result to higher dimensions.
\Remark{
\label{remark}
\begin{enumerate}
\item 
The proof will show that in each conformal class the set of metrics for 
which all non-zero eigenvalues are simple is residual.
\item 
The proof will also show that in any dimension the set of metrics for 
which $\dim E_{\lambda} \leq \operatorname{rank} \Sigma M$ for all non-zero 
$\lambda$ is residual. Only in dimension three is 
$\operatorname{rank} \Sigma M = 1$. 
\end{enumerate}
}

\section{Conformal deformations and Dirac eigenvalues}

Let $M$ be an $n$-dimensional manifold with a Riemannian metric $g$ and
a spin structure $\Spin(M,g) \rightarrow SO(M,g)$. We are going to study 
the deformation of $g$ given by $g^t = e^{2tf} g$, where $t$ is a real 
parameter and $f$ is a function on $M$.

There is an isomorphism $\gamma_{tf}: SO(M,g) \rightarrow SO(M,g^t)$ 
given by 
$$
\gamma_{tf}( \{e_1, \dots, e_n \} ) = \{e^{-tf} e_1, \dots, e^{-tf} e_n \}.
$$
Composition with $\gamma_{tf}$ turns the spin structure
$\Spin(M,g) \rightarrow SO(M,g)$ for $g$ into a spin structure for $g^t$;
$$
\Spin(M,g) \longrightarrow SO(M,g) \stackrel{\gamma_{tf}}{\longrightarrow} 
SO(M,g^t),
$$
and the associated spinor bundle 
$\Sigma M = \Spin(M,g) \times_{\rho} \Sigma$ is independent of $t$ and 
$f$. Let $D^t$ denote the Dirac operator acting on sections of 
$\Sigma M$ defined using the metric $g^t$. Then $D^t$ is related to 
$D= D^0$ by
$$
D^t \phi = e^{-tf} \left( 
D \phi + \frac{n-1}{2} t \grad f \cdot \phi
\right),
$$
where the Clifford multiplication by the gradient of $f$ is defined 
using the metric $g$, see for instance section 3.2.4. in \cite{baum81a}.
Define the Hilbert space $L^2(\Sigma M,g^t)$ as the completion of the 
smooth sections of $\Sigma M$ with respect to the inner product
$$
(\phi,\psi)_{g^t} = \int_M \la \phi, \psi \ra d \mu_{g^t}
= \int_M \la \phi, \psi \ra e^{ntf} d \mu.
$$
The proof of the following theorem is 
analogous to that of Theorem A.3 of \cite{bando-urakawa83a}, see also 
\cite{kato76a} and \cite{bourguignon-gauduchon92a}. 
\begin{thm} \label{perturbation-of-eigenspinors}
Suppose $\lambda$ is an eigenvalue of $D$ with $\dim_{\Co} E_{\lambda} = p$. 
Then there exist real analytic functions $\lambda_1^t, \dots, \lambda_p^t$ 
and curves of smooth spinor fields $\phi_1^t, \ldots, \phi_p^t$ 
(defined for small $t$) such that 
\begin{itemize}
\item   $D^t \phi_i^t = \lambda_i^t \phi_i^t$, $i=1,\ldots, p$,
\item   $\lambda_i^0 = \lambda$, $i=1,\ldots, p$, 
\item   $\phi_1^t, \ldots, \phi_p^t$ are real analytic as maps 
        $t \mapsto \phi_i^t \in L^2(\Sigma M,g)$,
\item   $\phi_1^t, \ldots, \phi_p^t$ are orthonormal in 
        $L^2(\Sigma M,g^t)$.
\end{itemize}
\end{thm}
\Remark{
If the spinor representation is quaternionic $\dim_{\Co}$ can be
replaced by $\dim_{\Qu}$, in which case ``orthonormal'' should be 
interpreted as
$$
(\phi_i,\phi_i)_g = 1, \quad (\phi_i,J \phi_i)_g =
(\phi_i,\phi_j)_g = (\phi_i,J \phi_j)_g =0, \quad i \neq j.
$$
}

Suppose that $t \mapsto (\lambda^t, \phi^t)$ is a one-parameter family 
of eigenspinors as given by Theorem~\ref{perturbation-of-eigenspinors};
$D^t \phi^t = \lambda^t \phi^t$ and $\| \phi^t \|_{g^t} = 1$ for small 
$t$. Then
\begin{eqnarray*}
\lambda^t &=& (D^t \phi^t, \phi^t )_{g^t} \\
&=&
\int_M \la e^{-tf} 
\left( D \phi^t + \frac{n-1}{2} t \grad f \cdot \phi^t \right),
\phi^t \ra e^{ntf} d\mu \\
&=&
\int_M \left( \la D \phi^t, \phi^t \ra e^{(n-1)tf}
+ \frac{n-1}{2} t \la \grad f \cdot \phi^t, \phi^t \ra e^{(n-1)tf}
\right) d\mu \\
&=&
\repart \int_M \la D \phi^t, \phi^t \ra e^{(n-1)tf} d\mu
\end{eqnarray*}
where the last equality follows since the second term is purely 
imaginary. Set $\lambda = \lambda^0$, $\phi = \phi^0$ and denote 
by prime the derivative with respect to $t$ at $t=0$. We have
\begin{eqnarray*}
\lambda' &=& 
\repart \int_M  \la D \phi', \phi \ra + \la D \phi, \phi' \ra
+ \la D \phi, \phi \ra (n-1) f d\mu \\
&=&
\lambda \int_M  2 \repart\la \phi', \phi \ra + (n-1)f|\phi|^2 d\mu
\end{eqnarray*}
and differentiating $\| \phi^t \|^2_{g^t} = 1$ at $t=0$ we get
$$
\int_M 2\repart \la \phi', \phi \ra + nf |\phi|^2 d\mu = 0
$$
which together give
\begin{equation} \label{lambda-prime}
\lambda' = - \lambda \int_M f|\phi|^2 d\mu.
\end{equation}

\section{Proof of Theorem \ref{Mainthm}}

Let $(M,g)$ be a three dimensional compact Riemannian spin manifold. 
The main technical point in the proof of Theorem \ref{Mainthm} is the 
following lemma.
\begin{lemma} \label{lemma-split}
Let $\lambda$ be a non-zero eigenvalue of $D_g$ with 
$\dim_{\Qu} E_{\lambda} = p > 1$. Then there is a conformal deformation of 
$g$ for which $E_{\lambda}$ splits into lower-dimensional eigenspaces. 
\end{lemma}
\begin{proof}
For each conformal deformation $g^t = e^{2tf} g$ 
Theorem~\ref{perturbation-of-eigenspinors} provides real analytic 
parametrizations $\lambda_1^t, \dots, \lambda_p^t$ of eigenvalues such 
that $\lambda_i^0 = \lambda$, $i=1,\ldots, p$. If for some $f$, $i$ and $j$ 
we have $\lambda_i^t \neq \lambda_j^t$ for all $t \neq 0$ we are done. 
Assume for a contradiction that this is not the case, but instead that 
$\lambda_i^t$ is independent of $i$, $t$ and $f$. Then for a given $f$ we 
can replace the eigenspinors $\phi_i^t$ from 
Theorem~\ref{perturbation-of-eigenspinors} by eigenspinors 
$\overline{\phi}_i^t = \sum_{j=1}^p U_{ij} \phi_j^t$ where $U_{ij}$ is a 
(constant) unitary matrix, which have the same properties. We may thus 
assume that $\phi_1^0$ and $\phi_2^0$ are the same for all conformal 
deformations.
 
For fixed $f$ and $p,q = 0,1$ let
$$
\alpha_{p,q}^t = 2^{-1/2}(\phi_1^t + i^p J^q \phi_2^t), 
\quad 
\beta_{p,q}^t = 2^{-1/2}(\phi_1^t - i^p J^q \phi_2^t). 
$$ 
Then 
$$
D \alpha_{p,q}^t = \lambda^t \alpha_{p,q}^t, \quad
D \beta_{p,q}^t = \lambda^t \beta_{p,q}^t,
$$
where $\lambda^t = \lambda_1^t = \lambda_2^t$. By (\ref{lambda-prime}) 
we have 
$$
\int_M f |\alpha_{p,q}^0|^2 d\mu = \lambda' = 
\int_M f |\beta_{p,q}^0|^2 d\mu.
$$
Since by assumption this holds for all $f$ and $\alpha_{p,q}^0$ and
$\beta_{p,q}^0$ are indedependent of $f$ we conclude that 
$|\alpha_{p,q}^0|^2 = |\beta_{p,q}^0|^2$ at each point. It follows that
$$
|\phi_1^0|^2 + |i^p J^q \phi_2^0|^2 
- 2 \repart \la \phi_1^0, i^p J^q \phi_2^0 \ra
= 
|\phi_1^0|^2 + |i^p J^q \phi_2^0|^2 
+ 2 \repart \la \phi_1^0, i^p J^q \phi_2^0 \ra
$$
so $\repart \la \phi_1^0, i^p J^q \phi_2^0 \ra = 0$ for $p,q = 0,1$, and
thus $\la \phi_1^0,\phi_2^0 \ra = \la \phi_1^0,J \phi_2^0 \ra = 0$ at 
each point. Since $\operatorname{rank}_{\Qu} \Sigma M = 1$ one of 
$\phi_1^0$ and $\phi_2^0$ has to vanish on an open set, and by unique 
continuation vanish identically. This is a contradiction which proves 
the lemma.
\end{proof}

Enumerate the non-zero eigenvalues of $D_g$ as 
$$
\ldots \leq \lambda_{-3} \leq \lambda_{-2} \leq \lambda_{-1} < 0 <
\lambda_{1} \leq \lambda_{2} \leq \lambda_{3} \leq \ldots
$$
with repetitions according to quaternionic multiplicity. Let $S^k(M)$, 
$k \geq 1$, be the subset of $R(M)$ consisting of Riemannian metrics 
for which $\lambda_{\pm 1}, \lambda_{\pm 2}, \dots, \lambda_{\pm k}$ 
are all different. \footnote{Comment added July 2013: 
For Proposition~\ref{prop-open-and-dense} and its proof to hold one must 
add the condition that the metrics in $S^k(M)$ have no harmonic spinors. 
This assumption is formulated slightly later here, in the proof of 
Theorem~\ref{Mainthm} on the next page. The fact that there is a problem 
in the published version of this paper was pointed out by Y.~Canzani and 
R.~Ponge, see the preprint \url{http://arxiv.org/abs/1207.0648}. 
The author wishes to express his thanks to them for their interest
and kind correspondence concerning this paper. Note that the argument 
works exactly as written to prove Remark~\ref{remark}; for perturbations 
within a conformal class the dimension of the kernel of $D$ is constant.}
\begin{prop} \label{prop-open-and-dense}
The sets $S^k(M)$, $k \geq 1$, are open and dense in $R(M)$.
\end{prop}
\begin{proof}
The eigenvalues of $D_g$ depend continuously on $g$, see Proposition~7.1 
in \cite{baer96b}. It follows that the $S^k(M)$ are open.

We prove that the $S^k(M)$ are dense by induction, for $k=1$ there is 
nothing to prove. Assume that $S^k(M)$ is dense in $R(M)$ and let $U$ be an 
open set in $R(M)$. We need to show that $U \cap S^{k+1}(M)$ is non-empty. 
By assumption we can find a metric $g \in U \cap S^k(M)$. Consider first 
the positive eigenvalues of $D_g$. The first $k$ are distinct, denote the 
multiplicity of $\lambda_k$ by $p$ so that 
$$
0 < \lambda_1 < \ldots < \lambda_{k-1} < \lambda_k
= \lambda_{k+1} = \ldots = \lambda_{k+p-1} < \lambda_{k+p} \leq \ldots .
$$
If $p=1$ we are done with the first step, so assume $p > 1$. Then by 
Lemma \ref{lemma-split} there is a conformal deformation which decreases 
the multiplicity of $\lambda_k$. By choosing the deformation parameter 
$t$ small enough we can guarantee to get a metric in $U \cap S^k(M)$ 
for which $\dim_{\Qu} E_{\lambda_k} < p$. Repeating this process we end
up with a metric $\hat{g} \in U \cap S^k(M)$ for which $\lambda_k$ has 
multiplicity $1$. 

The second step is to go through the same procedure with the negative 
eigenvalues of $D_{\hat{g}}$. When we do this we get a metric $\check{g}$ 
in $U \cap S^k(M)$ for which both $\lambda_k$ and $\lambda_{-k}$ have 
multiplicity $1$. The metric $\check{g}$ is thus in $U \cap S^{k+1}(M)$ 
and we have proved that $S^{k+1}(M)$ is dense in $R(M)$.  
\end{proof}

\begin{proof}[Proof of Theorem \ref{Mainthm}]
Let $S^0(M)$ denote the set of metrics on $M$ for which $\ker{D} = 0$.
Then $S(M)$ contains the intersection $\cap_{i=0}^{\infty} S^i(M)$. By 
Theorem~1.2 in \cite{maier97a} we know that $S^0(M)$ is open and dense 
in $R(M)$ and together with Proposition \ref{prop-open-and-dense} we 
conclude that $\cap_{i=0}^{\infty} S^i(M)$ is a residual set.
\end{proof}

\begin{proof}[Proof of Remark \ref{remark}]
We know from Lemma \ref{lemma-split} that we can split multiple 
eigenvalues within a conformal class. It follows that the intersection 
of $S^k(M)$ with a given conformal class is open and dense in the 
conformal class. This proves the first statement. For the second statement 
note that with the assumption 
$\dim E_{\lambda} > \operatorname{rank} \Sigma M$ the conclusion of 
Lemma~\ref{lemma-split} holds in any dimension. 
\end{proof}

\bibliographystyle{amsplain}

\providecommand{\bysame}{\leavevmode\hbox to3em{\hrulefill}\thinspace}

\end{document}